\documentclass[11pt,reqno]{amsart}

\usepackage{latexsym}
\usepackage{amsfonts}
\usepackage{amsmath}
\usepackage{url}
\usepackage{amssymb,amsmath}

\def\bl{\begin{lemma}}
\def\el{\end{lemma}}
\def\bth{\begin{theorem}}
\def\eth{\end{theorem}}
\def\bc{\begin{corollary}}
\def\ec{\end{corollary}}
\def\bcj{\begin{conjecture}}
\def\ecj{\end{conjecture}}
\def\bpr{\begin{proposition}}
\def\epr{\end{proposition}}
\def\bde{\begin{definition}}
\def\ede{\end{definition}}

\def\H{\mathbb{H}}

\def\QED{\hfill\qedsymbol}

\newcommand{\be}{\begin{eqnarray}}
\newcommand{\ee}{\end{eqnarray}}

\newcommand{\R}{{\mathbb R}}

\newcommand{\Z}{{\mathbb Z}}
\newcommand{\N}{{\mathbb N}}

\renewcommand{\and}{\hbox{ {\rm and} }}

\newtheorem{theorem}{Theorem}[section]
\newtheorem{definition}{Definition}[section]
\newtheorem{lemma}[theorem]{Lemma}

\newtheorem{corollary}[theorem]{Corollary}
\newtheorem{proposition}[theorem]{Proposition}
\newtheorem{conjecture}[theorem]{Conjecture}

\theoremstyle{definition}
\numberwithin{equation}{section}

\begin{document}
\title{Tightness of Fluctuations of
First Passage Percolation on Some Large Graphs}
\author{Itai Benjamini and Ofer Zeitouni}
\date{November 4, 2010}

\begin{abstract}
The theorem of Dekking and Host \cite{DH} regarding tightness
around the mean of  first passage percolation on the binary tree,
from the root to a boundary of a ball,
is generalized to a class of graphs which includes all lattices in hyperbolic spaces
and the lamplighter graph over $\N$.
This class of graphs is closed under product with
any bounded degree graph.
Few open problems and conjectures are gathered at the end.
\end{abstract}

\maketitle

\section{Introduction}

In {\it First Passage Percolation} (FPP) random i.i.d lengths are
assigned to the edges of a fixed graph. Among other questions one studies
the distribution of the distance from a fixed vertex to another vertex or
to a set, such as the boundary of a ball in the graph, see e.g. \cite{K}
for background.
Formally, given a rooted, undirected graph $G=(V,E)$ with root $o$,
let $D_n$ denote the collection of vertices at (graph) distance
$n$ from the root.
For $v\in D_n$, let ${\mathcal P}_v$ denote the collection
of
paths $(v_0=o,v_1,v_2,v_3,\ldots,v_k=v)$ (with $(v_{i-1},v_i)\in E$)
from $o$ to $v$.
Given a collection of positive i.i.d.
$\{X_e\}_{e\in E}$,
define, for $v\in E$,
\begin{equation}
	\label{eq-of1}
	Z_v=\min_{p\in {\mathcal P}_v} \sum_{e\in p} X_e\,.
\end{equation}
Because of the positivity assumption on the weights,
we may and will assume that any path in ${\mathcal P}_v$
visits each vertex of $G$
at most once.

For $n$ integer, let $Z_n^*=\min_{v\in D_n} Z_v$.
Under a mild moment condition on
the law of the random lengths $X_e$,
Dekking and Host \cite{DH} proved that for any regular tree,
$Z_n^*-EZ_n^*$,
the random distance from the root
to $D_n$ minus its mean,
is tight.
(Recall that a sequence of real valued random variables $\{X_n\}_{n\geq 0}$
is tight
iff for any $\epsilon > 0$, there is some $r_{\epsilon} \in \R$,
so that for all $n$, $P( |X_n|> r_{\epsilon}) < \epsilon$.)

We formulate here
a simple and general property of the underling graph $G$
and prove that for graphs satisfying
this property and a mild condition on
the law of $X_e$, the collection $\{Z_n^*-EZ_n^*\}_{n\geq 0}$
is tight.
Lattices in real hyperbolic spaces
$\H^d$,
the graph of the lamplighter over $\N$,
as well as
graphs of the form $G \times H$ where $G$
satisfies the conditions we list below
and $H$ is any bounded degree graph, are shown to possess
this property. (In passing, we mention that
the Euclidean case is wide open; it is known that
in two dimensions the fluctuations of
the distance are not tight,
see  \cite{NP,PP},  however only  very poor upper bounds
are known \cite{BKS}. For a special solved variant see \cite{J}.)
\medskip

In the next section we formulate the geometric condition on
the graph and the assumption on the distribution of the
edge weights $\{X_e\}_{e\in V}$;
we then state the tightness result, Theorem
\ref{theo-main},
 which is proved in Section \ref{sec-proof}.
We conclude  with a few open problems.

\section{A recursive structure in graphs and tightness}

Throughout,
let $dist_G$ denote the graph distance in $G$.
The following are the properties of $G$ and
the law of $X_e$ alluded to above.
\begin{itemize}
\item[(1)]
$G$ contains two vertex-disjoint subgraphs $G_1, G_2$,
which are isomorphic to $G$.
\item[(2)]
There exists $K<\infty$ so that $EX_e<K$,
and
$$
dist_G (Root_G, Root_{G_1}) = dist_G (Root_G, Root_{G_2}).
$$
%
\end{itemize}
One can replace Property (2) by the following.
\begin{itemize}
\item[(3)]
$X_e<K$ a.s., and
every vertex at distance $n$ from the root is connected to at least one
vertex of distance $n+1$.
\end{itemize}

Properties (1) and (2) imply that the binary tree
embeds quasi-isometrically into $G$, thus $G$ has exponential growth.
Property (3) is  called having "no dead ends" in geometric group theory terminology.
\medskip



\begin{theorem}
	\label{theo-main}
Assume Property (1) and either Property (2) or Property (3).
%
Then the sequence $\{(Z_n-EZ_n)\}_{n\geq 1}$ is tight.
\end{theorem}


\medskip

Note that a hyperbolic lattice in $\H^d, d \geq 2$,
intersected with a half space, admits the graph part of
Properties (1) and (2) above
(and probably (3) as well but we don't see a general proof).
This is due to topological transitivity of the action on the space of
geodesics, i.e. pairs of point of the boundary. There exist elements $g$ in
the authomorphism group of the hyperbolic space
that map the half space into arbitrarily small open sets of the boundary
and elements of this group map the lattice orbit to itself.
Note also
that by the
Morse lemma of hyperbolic geometry (see, e.g., \cite{BYI} p. 175), if one assumes in addition that
$X_e\geq \delta>0$ a.s. then
a path with minimal FPP length will be within a bounded
distance from a hyperbolic geodesic
and will not wind around, thus tightness for half space for weights
that are bounded below by a uniform positive constant
implies tightness for the whole space.
(Recall also that the regular tree is a lattice in $\H^2$; see \cite{L} for
some nice pictures of other planar hyperbolic lattices.)
\medskip

An example satisfying Properties (1) and (2) is given by the semi group of the lamplighter
over $\N$. Recall the graph of the lamplighter
over $\N$:
a vertex corresponds to a
scenery of $0$'s and $1$'s over $\N$, with finitely many $1$'s
with  a position of a lamplighter in $\N$;
edges either change the bit at the position of the
lamplighter or move the lamplighter one step to the left or the right,
see, e.g., \cite{L}.
If we fix the left most bit and restrict the lamplighter to integers
strictly bigger than $1$,
we get the required $G_0$ and $G_1$.

\medskip

It easy to see that if $G$ satisfies the properties in the theorem
then $G  \times H$ will too. In particular the theorem
applies to $T \times T'$ for two regular trees.
Note also that if $G$ satisfies the property (1) in the theorem,
then the lamplighter over $G$ will admit it as well.

\section{Proof of Theorem \ref{theo-main}}
\label{sec-proof}
The proof is based on a modification of an argument in \cite{DH}; a
related modification was used in \cite{BDZ}.
Note first that,
by construction,
$$
\mbox{(a)            }        EZ_{n+1}\geq EZ_n,
$$
because to get to distance $n+1$ a path has
to pass through distance $n$ and
the weights $\{X_e\}$ are positive.

Under Property (3), one has  in addition

(a') $Z_n$ and $Z_{n+i}$ can be constructed on the same space so that
$$
    Z_{n+i}\geq Z_n           \mbox{      while        }  Z_{n+i}\leq Z_n+Ki.
$$
(The first inequality does not need Property
(3), but the second does --- one just goes
forward from the minimum at distance $n$, $i$ steps.)

On the other hand,
from Property (1),
$$
EZ_{n+1} \leq E(\min (Z_{n-R_1+1}, Z_{n-R_2+1}')+K C,
$$
where $R_i=dist_G(Root_G,Root_{G_i})$, $C=\max(R_1,R_2)$, and
$Z'_m$ denotes a identically distributed independent copy of $Z_m$ .

Since $\min(a, b) = {a+b \over 2} - {|a-b| \over 2}$,
$$
EZ_{n+1}\leq (1/2) [EZ_{n-R_1+1}+ EZ_{n-R_2+1}'-E|Z_{n-R_1+1}-Z_{n-R_2+1}'|]
+KC.
$$
Therefore, with $n_i=n+1-R_i$,
$$
E|Z_{n_1}-Z_{n_2}'|\leq [-2EZ_{n+1}+EZ_{n_1}+EZ_{n_2}]+2KC
.$$
If $R_1=R_2$ (i.e. Property (2) holds), then, using (a),
$$
E|Z_{n_1}-Z_{n_1}'|\leq 2KC\,,
$$
and the tightness follows by standard arguments (see, e.g., \cite{DH}).
Otherwise,
assume Property (3) with
$n_2>n_1$. By (a'), we can construct a version of $Z_{n_1}'$,
independent
of $Z_{n_1}$,
so that $|Z_{n_2}'-Z_{n_1}'|\leq K(n_2-n_1)$. Therefore,
$$
E|Z_{n_1}-Z_{n_1}'|\leq [-2EZ_{n+1}+EZ_{n_1}+EZ_{n_2}] +2KC+K(R_1-R_2).
$$
Applying again (a) we get, for some constant $C'$,
$$
E|Z_{n_1}-Z_{n_1}'|\leq 2KC+K(R_1-R_2)
\leq  C'K,
$$
and as before it is standard that this implies tightness.
\QED

\section{Questions}

\noindent
{\bf Question 1:}
Extend the theorem to the
lamplighter group over $\Gamma$, for any finitely generated group
$\Gamma$; start with $\Z$.

\noindent
{\bf Question 2:} Show that tightness of fluctuations
is a quasi-isometric invariant.
In particular, show this in the class of Cayley graphs.

\noindent
{\bf Question 3:} The lamplighter over $\Z$ is a rather small
group among the finitely generated groups with exponential growth.
It is solvable, amenable and Liouville.
This indicates that all Cayley graphs of exponential growth are tight.
We ask then which Cayley graphs admit tightness;
is there an infinite Cayley graph,
which is not quasi-isometric to $\Z$ or $\Z^2$, for which tightness
does not hold? Start with  a sub exponential example with tightness
or even only variance smaller than on $\Z^2$.

\noindent
{\bf Question 4:}
(Gabor Pete) Note that requiring (1) only quasi-isometrically (plus the root
condition of (2)) does not imply exponential growth, because
when one iterates,
one may collect a factor (from quasiness) each time, killing
the exponential growth.
E.g., branch groups like Grigorchuk's group \cite{GP}, where
$G \times G$ is a subgroup of $G$, may have intermediate
growth, see e.g. \cite{N}.
This condition is somewhat in the spirit of property (1).
Bound the variance for FPP on the Grigorchuk's group.

Maybe ideas related to the one above will be
useful in proving at least a sublinear variance?
\medskip

The last two questions are regarding point to point FPP.

\noindent
{\bf Question 5:} We conjecture that in any hyperbolic
lattice the point to point FPP
fluctuations admit a central
limit theorem with variance proportional to the distance.
This is motivated by the fact that,
due to the Morse lemma,
the minimal path will be in a bounded neighborhood of the hyperbolic
geodesic, and for cylinders a CLT is  known to hold \cite{CD}.
\medskip

A related question is the following.
Assume that for any pair of vertices in a Cayley graph the variance
of point to point FPP is proportional to the distance,
is the Cayley graph hyperbolic? Alternatively,
what point to point variances can be achieved
for Cayley graphs? As pointed out above,
the only behavior known is
linear in the distance (for $\Z$), the conjectured (and proved in some cases)
behavior
 for $\Z^2$,
which is
the distance to the power $2/3$. Can the bound  or proof of
theorem \ref{theo-main} be adapted to give point to point order $1$ variance
for $T \times \Z^d$ or $T \times T$ or some other graphs? Are other behaviors possible?

\noindent
{\bf Question 6:} In \cite{BHH}, among other things,
tightness was proved for point to point 
FPP  between random
vertices   in the configuration model of random $d$-regular graph with
exponential weights.
Does tightness hold for more general weights, or 
for point to point FPP  between random
vertices on expanders?
\medskip

All the questions above are regarding the second order issue of bounding fluctuations.
The fundamental fact regarding FPP on $\Z^d$ is the shape theorem, see e.g. \cite{K}.
That is, rescale the random FPP metric then the limiting metric space a.s. exists
and is $\R^d$ with some deterministic norm. The subadditive ergodic theorem is a key in the proof.
We conjecture that FPP on Cayley graph of groups of polynomial growth also admits a shape theorem.
What can replace the subadditive ergodic theorem in the proof? Start with

\noindent
{\bf Question 7:} Prove a shape theorem for FPP on the Cayley graph of the discrete Heisenberg group.

\medskip

\noindent {\bf Acknowledgements:} Thanks to Pierre Pansu and Gabor Pete for very useful discussions, and to R. van der Hofstad for bringing 
\cite{BHH} to our attention.

 \vspace{.01 in}\noindent
 The Weizmann Institute of Science, \\
 Rehovot  76100, Israel.

\end{document}